# USING THE BOOTSTRAP TO QUANTIFY THE AUTHORITY OF AN EMPIRICAL RANKING

BY PETER HALL AND HUGH MILLER

*University of Melbourne*

The bootstrap is a popular and convenient method for quantifying the authority of an empirical ordering of attributes, for example of a ranking of the performance of institutions or of the influence of genes on a response variable. In the first of these examples, the number, $p$, of quantities being ordered is sometimes only moderate in size; in the second it can be very large, often much greater than sample size. However, we show that in both types of problem the conventional bootstrap can produce inconsistency. Moreover, the standard $n$-out-of-$n$ bootstrap estimator of the distribution of an empirical rank may not converge in the usual sense; the estimator may converge in distribution, but not in probability. Nevertheless, in many cases the bootstrap correctly identifies the support of the asymptotic distribution of ranks. In some contemporary problems, bootstrap prediction intervals for ranks are particularly long, and in this context, we also quantify the accuracy of bootstrap methods, showing that the standard bootstrap gets the order of magnitude of the interval right, but not the constant multiplier of interval length. The $m$-out-of-$n$ bootstrap can improve performance and produce statistical consistency, but it requires empirical choice of $m$; we suggest a tuning solution to this problem. We show that in genomic examples, where it might be expected that the standard, "synchronous" bootstrap will successfully accommodate nonindependence of vector components, that approach can produce misleading results. An "independent component" bootstrap can overcome these difficulties, even in cases where components are not strictly independent.

**1. Introduction.** The ordering of a sequence of random variables is often a major aspect of contemporary statistical analyses. For example, data on the comparative performance of institutions (e.g., local governments, or health providers, or universities) are frequently summarized by reporting the ranking of empirical values of a performance measure; and the relative influence of genes on a particular response is sometimes indicated by ranking the values of the weights that are applied to them after the application of a variable selector, such as the lasso. It is reasonable to argue that, especially in contentious situations, no ranking should be unaccompanied by a measure of its authority [Goldstein and Spiegelhalter (1996)]. The bootstrap is a popular approach to developing such a measure.

In this paper, we report on both the theoretical and the numerical properties of bootstrap estimators of the distributions of rankings. We show that the standard $n$-out-of-$n$ bootstrap generally fails to give consistency, and in fact may not produce distribution estimators that converge either almost surely or in probability. The $m$-out-of-$n$ bootstrap overcomes these difficulties, but requires empirical choice of $m$. We suggest a tuning approach to solving this problem. This technique remains appropriate in cases where the number, $p$ say, of populations is very large, although in that context, one could also regard $m$ as a means of setting the level of sensitivity of the bootstrap to near-ties among ranks, rather than as a smoothing parameter.

In some contemporary prediction problems, the empirical rank is quite highly variable. We develop mathematical models in this setting, and explore the validity of bootstrap methods there. In particular, we show that the inherent inconsistency of the standard $n$-out-of-$n$ bootstrap does not prevent that method from correctly capturing the order of magnitude of the expected value of rank, or the expected length of prediction intervals, although it leads to errors in estimators of the constant multiplier of that order of magnitude.

Another issue is that of adequately reflecting, in the bootstrap algorithm, dependence among the datasets representing the different populations, for example, data on the performances of different health providers, or on the expression levels of different genes. In examples of the first type, where different institutions are being ranked, the assumption of independence is often appropriate; it can usually be accommodated through conditioning. In such cases, resampling can be implemented in a way that explicitly reflects population-wise independence.

However, in the genomic example, data on expression levels of different genes from the same individual are generally not independent. In this setting, using the standard nonparametric bootstrap to assess the authority of ranking would seem to be a good choice, since in more conventional problems it captures well the dependence structure of data vectors. However, we show that, even when the number of variables being ranked is much less than sample size, the standard approach can give unreliable results in some problems. This is largely because knowing the composition of a resample for the $j$th population (e.g., for the $j$th gene, in the genomic example) identifies exactly the resamples for other genes. Therefore, the resamples for different populations are hardly independent, even conditional on the original data.

This has a variety of repercussions. For example, it implies that standard bootstrap probabilities, when computed conditional on the information we have in the resample about the $j$th gene, degenerate to indicator functions. Conditional inference is attractive in ranking problems, since it can lead to substantial reductions in variability. To overcome the problem, we suggest using an "independent component" version of the bootstrap, where the bootstrap is applied as though the ranked variables were statistically independent. This approach can be valid even in the case of nonindependence. (In order to make it clear that in this setting we use the term "standard bootstrap" to mean the resampling of $p$-vectors of data, we shall refer to this bootstrap method as the "synchronous" bootstrap; the standard bootstrap results in vector components being synchronised with one another in each resampling step.)

It is possible to generalize our treatment to cases where several rankings are undertaken jointly, for example, where universities are ranked simultaneously in terms of the quality of their graduate programs and the career prospects of their undergraduates. Our main conclusions about the relative merits of different bootstrap methods persist in this more general setting, although a detailed treatment of that case would be significantly longer and more complex.

Work on the bootstrapping of statistics related to ranks includes that of Srivastava ([1987](#)), who introduced bootstrap methods for a class of ranking and slippage problems (although not directly related to the problems discussed in this paper); Tu, Burdick and Mitchell ([1992](#)), who discussed bootstrap methods for canonical correlation analysis; Larocque and Léger ([1994](#)), Steland ([1998](#)) and Peilin et al. ([2000](#)), who developed bootstrap methods for quantities such as rank tests and rank statistics; Goldstein and Spiegelhalter ([1996](#)), who discussed bootstrap methods for constructing interval estimates; Langford and Leyland ([1996](#)), who addressed bootstrap methods for ranking the performance of doctors; Cesário and Barreto ([2003](#)), Hui, Modarres and Zheng ([2005](#)) and Taconeli and Barreto ([2005](#)), who discussed bootstrap methods for ranked set sampling; and Mukherjee et al. ([2003](#)), who developed methods for gene ranking using bootstrapped $p$-values.

## 2. Model and methodology.

2.1. *Model.* Assume we have datasets $\mathcal{X}_1, \ldots, \mathcal{X}_p$ drawn from populations $\Pi_1, \ldots, \Pi_p$, respectively, and that for the $j$th population there is an associated parameter $\theta_j$ which measures, for example, the strength of an attribute in the population, or the performance of an individual or an organisation related to the population, or the esteem in which an institution or a program is held. If the $\theta_j$'s were known, then our ranking of the populations would be

$$(2.1) \qquad r_j = 1 + \sum_{k \neq j} I(\theta_k \geq \theta_j) \qquad \text{for } j = 1, \ldots, p,$$

say, signifying that $r_j$ is the rank of the $j$th population. Here, tied rankings can be considered to have been broken arbitrarily, for example at random.

We wish to develop an empirical version of the ranking at (2.1). For this purpose, we compute from $\mathcal{X}_j$ an estimator $\hat{\theta}_j$ of $\theta_j$, for $1 \leq j \leq p$, and we rank the populations in terms of the values of $\hat{\theta}_j$. In particular, if we have $\hat{\theta}_1, \ldots, \hat{\theta}_p$, then we write

$$(2.2) \qquad \hat{r}_j = 1 + \sum_{k \neq j} I(\hat{\theta}_k \geq \hat{\theta}_j) \qquad \text{for } j = 1, \ldots, p,$$

to indicate the empirical version of (2.1). Again, ties can be broken arbitrarily, although in the case of (2.2) the noise implicit in the estimators $\hat{\theta}_j$ often means that there are no exact ties.

We shall treat two cases: "fixed $p$" and "large $p$," distinguished in theoretical models by taking $p$ fixed and allowing $n$ to diverge, and by permitting $p$ to diverge, respectively. Cases covered by the latter model include instances where $\mathcal{X}_0$ is a set of $p$-vectors, say $\mathcal{X}_0 = \{X_1, \ldots, X_n\}$ where $X_i = (X_{i1}, \ldots, X_{ip})$. There, $\mathcal{X}_j = \{X_{1j}, \ldots, X_{nj}\}$ is the set of $j$th components of each data vector, and in particular each $\mathcal{X}_j$ is of the same size. This example arises frequently in contemporary problems in genomics, where $X_i$ is the vector of expression-level data on perhaps $p = 5000$ to 20,000 genes for the $i$th individual in a population. In such cases, $n$ can be relatively small, for example, between 20 and 200. The vectors $X_i$ can generally be regarded as independent, but not so the components $\mathcal{X}_1, \ldots, \mathcal{X}_p$. However, as we shall argue in Section 4.1, there may be advantages in conducting inference as though the components were independent, even when that assumption is incorrect.

2.2. *Basic bootstrap methodology.* The authority of the ranking at (2.2), as an approximation to that at (2.1), can be queried. A simple approach to quantifying the authority is to repeat the ranking many times in the context of bootstrap resamples $\mathcal{X}_1^*, \ldots, \mathcal{X}_p^*$, which replace the respective datasets $\mathcal{X}_1, \ldots, \mathcal{X}_p$. In particular, for each sequence $\mathcal{X}_1^*, \ldots, \mathcal{X}_p^*$, we can compute the respective versions $\hat{\theta}_1^*, \ldots, \hat{\theta}_p^*$ of the estimators of $\theta_j$, and calculate the bootstrap version of (2.2):

$$(2.3) \qquad \hat{r}_j^* = 1 + \sum_{k \neq j} I(\hat{\theta}_k^* \geq \hat{\theta}_j^*) \qquad \text{for } j = 1, \ldots, p.$$

The bootstrap here can be of conventional $n$-out-of-$n$ type, either parametric or nonparametric, or it can be the $m$-out-of-$n$ bootstrap (again either parametric or nonparametric), where the resamples $\mathcal{X}_j^*$ are of smaller size than the respective samples $\mathcal{X}_j$. For definiteness, in Section 4, where we need to refer explicitly to the implementation of bootstrap methods, we shall use the nonparametric bootstrap. However, our conclusions also apply to parametric bootstrap methods. More generally, the way in which the bootstrap resamples $\mathcal{X}_j^*$ are constructed can depend on the nature of the data. See Section 4.1 for discussion.

One question in which we are obviously interested is whether the bootstrap captures the distribution of $\hat{r}_j$ reasonably well, for example, whether

$$(2.4) \qquad P(\hat{r}_j^* \leq r | \mathcal{X}) - P(\hat{r}_j \leq r) \to 0,$$

in probability for each integer $r$, as $n \to \infty$. The answer to this question, if we use the familiar $n$-out-of-$n$ bootstrap, is generally "only in cases where the limiting distribution of $\hat{r}_j$ is degenerate." However, the answer is more positive if we employ the $m$-out-of-$n$ bootstrap. There, if the populations $\Pi_1, \ldots, \Pi_p$ are kept fixed in an asymptotic study then

$$(2.5) \qquad \begin{array}{l} \text{the limiting distribution of } \hat{r}_j \text{ is supported on the set of integers } \{k_1 + 1,\ k_1 + 2, \ldots, k_2\}, \text{ where } k_1 = \sum_k I(\theta_k > \theta_j) \text{ and } k_2 = \sum_k I(\theta_k \geq \theta_j), \end{array}$$

and the $m$-out-of-$n$ bootstrap consistently estimates this distribution. In particular, (2.4) holds; see Section 3 for details. However (still in the case of fixed $p$), if we are more ambitious and permit the population distributions to vary with $n$ in such a way that the limiting distribution is more complex than that prescribed by (2.5), then even the $m$-out-of-$n$ bootstrap may fail to give consistency.

Having computed a bootstrap approximation $P(\hat{r}_j^* \leq r | \mathcal{X})$ to the probability $P(\hat{r}_j \leq r)$, we can calculate an empirical approximation to a prediction interval, specifically an interval $[r_1, r_2]$ within which $\hat{r}_j$ lies with given probability, for example, 0.95. Goldstein and Spiegelhalter (1996) refer to such intervals as "overlap intervals," since they are generally displayed in a figure which shows the extent to which they overlap. Particularly, when $p$ is relatively small, the discrete nature of the distribution of $\hat{r}_j$ makes it a little awkward to discuss the accuracy of bootstrap prediction intervals, and so we focus instead on measures of the accuracy of distributional approximations, for example, (2.4) and (2.5).

**3. The case of $p$ distinct populations.**

3.1. *Preliminary discussion.* Write $n_j$ for the size of the sample $\mathcal{X}_j$. The values of $n_j$ may differ, but we shall assume that they all of the same order. That is, writing $n = p^{-1} \sum_j n_j$ for the average sample size, we have:

$$(3.1) \qquad n^{-1} \sup_{1 \le j \le p} n_j = O(1), \qquad 1 = O\Big(n^{-1} \inf_{1 \le j \le p} n_j\Big).$$

When interpreting (3.1) it is convenient to think of $n$ as the "asymptotic parameter," for example, the quantity which we take to diverge to infinity, and to consider $n_1, \ldots, n_p$ as functions of $n$.

When using the $m$-out-of-$n$ bootstrap, where a resample of size $m_j < n_j$ is drawn either from the population distribution with estimated parameters (the parametric case) or by with-replacement resampling from the sample $\mathcal{X}_j$ (the case of the nonparametric bootstrap), and $\mathcal{X}_j$ is of size $n_j$, we assume that the average resample size, $m = p^{-1} \sum_j m_j$, satisfies the analogue of (3.1):

$$(3.2) \qquad m^{-1} \sup_{1 \le j \le p} m_j = O(1), \qquad 1 = O\Big(m^{-1} \inf_{1 \le j \le p} m_j\Big).$$

Furthermore, we ask that $m$ be large but $m/n$ be small.

In the cases of both fixed and divergent $p$, the properties of $\hat{r}_j$ and $\hat{r}_j^*$ are strongly influenced by the potential presence of tied values of $\theta_j$. However, it is perhaps unreasonable to assume, in practice, that two values of $\theta_j$ are exactly tied, although there might be cases where two values are so close that, for most practical purposes, the properties of $\hat{r}_j$ for small to moderate $n$ are similar to those that would occur if the values were tied. The borderline case is that where two values of $\theta_j$ differ by only a constant multiple of $n^{-1/2}$, with $n$ denoting average sample size. (This requires the distribution of the populations $\Pi_j$ to vary with $n$.) If the constant is sufficiently large, then practically speaking, the two values of $\theta_j$ are not tied, but if the constant is small then a tie might appear to be present.

To reflect this viewpoint, we shall for any particular $j$ and for all $k \ne j$, write

$$(3.3) \qquad \theta_k = \theta_j + n^{-1/2} \omega_{jk},$$

where the $\omega_{jk}$'s are permitted to depend on $n$. Of course, (3.3) amounts to a definition of $\omega_{jk}$, and if the quantities $\theta_k$, for $1 \le k \le p$, are all fixed then (3.3) implies that $\omega_{jk}$ either vanishes or diverges to either $+\infty$ or $-\infty$, in the latter two cases in proportion to $n^{1/2}$. However, since we shall permit the distributions of the populations $\Pi_k$, and hence also the $\theta_k$'s, to depend on $n$, then the problem can be set up in such a way that the $\omega_{jk}$'s have many different modes of behavior.

In the case of the $m$-out-of-$n$ bootstrap, where $m \to \infty$ but $m/n \to 0$, sensitivity is somewhat reduced by using a smaller resample size. Reflecting this restriction, in the $m$-out-of-$n$ bootstrap setting, we use the following formula to define quantities $\omega'_{jk}$, in place of the $\omega_{jk}$'s at (3.3):

$$(3.4) \qquad \theta_k = \theta_j + m^{-1/2}\omega'_{jk}.$$

It can be proved that, under regularity conditions, the sum over $r$ of the squared distance between the $m$-out-of-$n$ bootstrap approximation to the distribution function of $\hat{r}_j$, and the limiting form $G_j$ of that distribution [see (3.9) below], equals $C_1 m^{-1} + C_2 mn^{-1} + o(m^{-1} + mn^{-1})$, where $C_1$ and $C_2$ are positive constants. This result implies that the asymptotically optimal choice of $m$ equals the integer part of $(C_1 n/C_2)^{1/2}$. However, this limit-theoretic argument is not always valid when $p$ is large, and even in the case of small $p$ it is not straightforward to estimate the ratio $C_1/C_2$. In Section 3.4, we suggest an alternative, relatively flexible, method for choosing $m$.

In most cases, where there are $p$ distinct populations, it is reasonable to argue that the datasets $\mathcal{X}_1, \ldots, \mathcal{X}_p$ are independent. For example, $\mathcal{X}_j$ might represent a sample relating to the performance of the $j$th of $p$ health providers that are being operated essentially independently [see, e.g., Goldstein and Spiegelhalter (1996)], and the data in $\mathcal{X}_j$ would be gathered in a way that is largely independent of data for other health providers. To the extent to which the data are related, for example, through the common effects of government policies, or shared health-care challenges such as epidemics, we might interpret our analysis as conditional on those effects.

If the assumption of independence is valid, then it is straightforward to reflect the assumption during the resampling operation, obtaining bootstrap parameter estimators $\hat{\theta}_1^*, \ldots, \hat{\theta}_p^*$ that are independent conditional on $\mathcal{X} = \bigcup_j \mathcal{X}_j$. If the independence assumption is not appropriate, then resampling is generally a more complex operation, and may be so challenging as to be impractical. In the remainder of this section, we shall assume that $\mathcal{X}_1, \ldots, \mathcal{X}_p$ are independent, and that $\hat{\theta}_1^*, \ldots, \hat{\theta}_p^*$ are independent conditional on $\mathcal{X}$.

Sections 3.2 and 3.5 will outline theoretical properties in the case of fixed $p$ and increasingly large $p$, respectively. To simplify and abbreviate our discussion we shall state our main results only for one $j$ at a time, but joint distribution properties can also be derived, analogous to those in Theorem 4.1.

### 3.2. *Theoretical properties in the case of fixed $p$.*

To set the scene for our results, we note first that, under mild regularity conditions, it holds true that for fixed $p$, for each $1 \leq j \leq p$ and for each real number $x$,

$$(3.5) \qquad \begin{aligned} P\{n^{1/2}(\hat{\theta}_j - \theta_j) \leq \sigma_j x\} &\to \Phi(x), \\ P\{m^{1/2}(\hat{\theta}_j^* - \hat{\theta}_j) \leq \sigma_j x | \mathcal{X}\} &\to \Phi(x), \end{aligned}$$

where the asymptotic standard deviations $\sigma_j \in (0, \infty)$ do not depend on $n$, $\Phi$ denotes the standard normal distribution function, and the convergence in the second part of (3.5) is in probability. In that second part, the value of $m$ equals $n$ if we are using the conventional bootstrap, and equals $m$ if we are using the $m$-out-of-$n$ bootstrap.

The first formula in (3.5) is the conventional statement that the statistics $\hat{\theta}_j$ are asymptotically normally distributed, and the second is the standard bootstrap form of that assumption. It asserts only that the bootstrap estimator of the distribution of $n^{1/2}(\hat{\theta}_j - \theta_j)$ is consistent for the normal distribution with zero mean and variance $\sigma_j^2$.

In this section, we keep $p$ fixed as we vary $n$, although we permit the distributions of the populations $\Pi_1, \ldots, \Pi_p$ to depend on $n$. Let $N_1, \ldots, N_p$ denote independent standard normal random variables and, given constants $c_1, \ldots, c_p$, let $F_j(\cdot | c_1, \ldots, c_p)$ denote the distribution function of the random variables

$$1 + \sum_{k \,:\, k \neq j} I(\sigma_j N_j \leq \sigma_k N_k + c_k).$$

The value of $c_j$ has no influence on $F_j$, but it is cumbersome to reflect this in notation.

THEOREM 3.1.  *Assume that $p$ is fixed and the datasets $\mathcal{X}_1, \ldots, \mathcal{X}_p$ are independent, that $\hat{\theta}_1^*, \ldots, \hat{\theta}_p^*$ are independent conditional on $\mathcal{X}$, and that (3.1), (3.2) (if using the $m$-out-of-$n$ bootstrap) and (3.5) hold. [In (3.5), we take $m = n$ unless using the $m$-out-of-$n$ bootstrap.]*

(i) *For each integer $r$,*

$$(3.6) \qquad\qquad P(\hat{r}_j \leq r) - F_j(r | \omega_{j1}, \ldots, \omega_{jp}) \to 0$$

*as $n \to \infty$.*

(ii) *Using the standard $n$-out-of-$n$ bootstrap, either parametric or nonparametric, define the $\omega_{jk}$'s by (3.3). Then there exists a sequence of random variables $Z_1, \ldots, Z_p$, depending on $n$ and being, for each choice of $n$, independent and having the standard normal distribution, such that*

$$(3.7) \quad P(\hat{r}_j^* \leq r | \mathcal{X}) - F_j(r | \omega_{j1} + \sigma_1 Z_1 - \sigma_j Z_j, \ldots, \omega_{jp} + \sigma_p Z_p - \sigma_j Z_j) \to 0$$

*in probability as $n \to \infty$.*

(iii) *In the case of the $m$-out-of-$n$ bootstrap, again either parametric or nonparametric, and for which $m/n \to 0$ and $m \to \infty$, define the $\omega_{jk}'$'s by (3.4). Then (3.7) alters to*

$$(3.8) \qquad\qquad P(\hat{r}_j^* \leq r | \mathcal{X}) - F_j(r | \omega_{j1}', \ldots, \omega_{jp}') \to 0$$

*in probability as $n \to \infty$.*

3.3. *Interpretation of Theorem 3.1.* To illustrate the implications of the theorem, let us assume that $\omega_{jk}$, defined by (3.3), has (for each $j$ and $k$) a well-defined limit (either finite or infinite) as $n \to \infty$, and that $\omega_{jk} \to +\infty$ for $k \in \mathcal{K}_+$, $\omega_{jk} \to -\infty$ for $k \in \mathcal{K}_-$, and $\omega_{jk}$ has a finite limit, $\omega_{jk}^0$ say, for $k \in \mathcal{K}_j = \{1, \ldots, p\} \setminus (\{j\} \cup \mathcal{K}_+ \cup \mathcal{K}_-)$. (Both $\mathcal{K}_+$ and $\mathcal{K}_-$ may depend on $j$.) Define $G_j$ to be the distribution function of

$$1 + (\#\mathcal{K}_+) + \sum_{k \in \mathcal{K}_j} I(\sigma_j N_j \leq \sigma_k N_k + \omega_{jk}^0).$$

Then $F_j(r|\omega_{j1}, \ldots, \omega_{jp}) \to G_j(r)$, and so (3.6) implies that, as $n \to \infty$,

$$(3.9) \qquad P(\hat{r}_j \leq r) \to G_j(r)$$

for each integer $r$.

Analogously to the argument leading from (3.6) to (3.9), result (3.7) implies that, in the case of the $n$-out-of-$n$ bootstrap,

$P(\hat{r}_j^* \leq r | \mathcal{X})$ converges in distribution to the random variable,

$$(3.10) \quad P\bigg[1 + (\#\mathcal{K}_+) + \sum_{k \in \mathcal{K}_j} I\{\sigma_j(N_j + N_j') \leq \sigma_k(N_k + N_k') + \omega_{jk}^0\}$$

$$\leq r | N_1, \ldots, N_p\bigg],$$

where $N_1, \ldots, N_p, N_1', \ldots, N_p'$ are independent standard normal random variables. However, the convergence of $P(\hat{r}_j^* \leq r | \mathcal{X})$ is not in probability.

If $\mathcal{K}_+ \cup \mathcal{K}_- = \{1, \ldots, p\} \setminus \{j\}$, which occurs for example, if the $\theta_k$'s are fixed and there are no ties for the value of $\theta_j$, then it follows from (3.6) and (3.7) that $P(\hat{r}_j = r_j) \to 1$ and $P(\hat{r}_j^* = r_j | \mathcal{X}) \to 1$ in probability, where $r_j$ denotes the rank of $\theta_j$ in the set of all $\theta_k$'s. Therefore, in this degenerate setting, the standard $n$-out-of-$n$ bootstrap correctly captures the asymptotic distribution of $\hat{r}_j$.

In all other cases, however, the limiting distribution of $\hat{r}_j$ [see (3.9)] does not equal the limit of the $n$-out-of-$n$ bootstrap distribution of $\hat{r}_j^*$ [see (3.10)]. Nevertheless, it is clear from (3.9) and (3.10) that

(3.11) The support of the limiting distribution of $\hat{r}_j$, and the support of the weak limit of the distribution of $\hat{r}_j^*$ given $\mathcal{X}$, are identical, and both are equal to the set $\{\#\mathcal{K}_+ + 1, \ldots, \#\mathcal{K}_+ + \#\mathcal{K}_j + 1\}$.

To this extent the standard $n$-out-of-$n$ bootstrap correctly captures important aspects of the distribution of $\hat{r}_j$.

Superficially, (3.8) seems to imply that the $m$-out-of-$n$ bootstrap overcomes this problem. However, the $\omega_{jk}'$'s are now defined by (3.4), and are different from the $\omega_{jk}$'s at (3.3). As a result, the $m$-out-of-$n$ bootstrap does not, in general, correctly capture the limiting distribution at (3.9). Nevertheless, if

$$\text{(3.12)} \qquad \text{for each } k \neq j \qquad \text{either } m^{1/2}(\theta_k - \theta_j) \to \pm\infty \quad \text{or}$$
$$n^{1/2}(\theta_k - \theta_j) \to 0,$$

then $P(\hat{r}_j^* \leq r | \mathcal{X}) - P(\hat{r}_j \leq r) \to 0$ in probability, that is, (2.4) holds. In particular, the $m$-out-of-$n$ bootstrap consistently estimates the distribution of empirical ranks. Under condition (3.12), the following analogue of (3.11) holds for the $m$-out-of-$n$ bootstrap:

(3.13)    The limiting distributions of $\hat{r}_j$, and of $\hat{r}_j^*$ conditional on $\mathcal{X}$, are identical when using the $m$-out-of-$n$ bootstrap, and the support of each equals the set $\{\#\mathcal{K}_+ + 1, \ldots, \#\mathcal{K}_+ + \#\mathcal{K}_j + 1\}$.

Property (3.12) holds if the $\theta_k$'s are all fixed (i.e., do not depend on $n$). Therefore, the $m$-out-of-$n$ bootstrap correctly estimates the distribution of ranks in the presence of ties, when the populations are kept fixed as sample sizes diverge, and also in other cases where the differences $\theta_k - \theta_j$ are of either strictly larger order than $m^{-1/2}$ or strictly smaller order than $n^{-1/2}$. When (3.12) holds, the asymptotic distribution of $\hat{r}_j$ is supported on a set the size of $\#\mathcal{K}_j$, that is the number of integers $k$ for which $m^{1/2}(\theta_k - \theta_j) \to 0$.

3.4. *Methods for choosing $m$.* Consider a comparison of two of the populations $\Pi_j$ and $\Pi_k$, and focus on the probability of ranking one higher than the other using the $m$-out-of-$n$ bootstrap. Assuming (3.5) and letting $c = (\sigma_j^2 + \sigma_k^2)^{-1/2}$, we see that

$$P(\hat{r}_j^* < \hat{r}_k^* | \mathcal{X}) - P(\hat{r}_j < \hat{r}_k)$$
$$= P(\hat{\theta}_j^* > \hat{\theta}_k^* | \mathcal{X}) - P(\hat{\theta}_j > \hat{\theta}_k)$$
$$\approx \Phi\{m^{1/2}c(\hat{\theta}_j - \hat{\theta}_k)\} - \Phi\{n^{1/2}c(\theta_j - \theta_k)\}$$
$$\approx \Phi\{m^{1/2}c(\theta_j - \theta_k) + c(m/n)^{1/2}Z\} - \Phi\{n^{1/2}c(\theta_j - \theta_k)\}$$
$$= \Phi\{(m/n)^{1/2}(-c\omega_{jk} + Z)\} - \Phi(-c\omega_{jk}).$$

Here, $Z$ denotes a realization of a normal random variable, and $\Phi$ is the standard normal distribution function. Thus, choosing $m$ to minimize the squared difference between the bootstrapped and true probabilities is approximately equivalent to choosing $m$ to minimize the expression

$$\text{(3.14)} \qquad [\Phi\{(m/n)^{1/2}(-c\omega_{jk} + Z)\} - \Phi(-c\omega_{jk})]^2.$$

If $\omega_{jk} \to \pm\infty$, then the expression is minimized as long as $(m/n)^{1/2}\omega_{jk} \to \pm\infty$ too, which guarantees that $m \to \infty$ as long as $\omega_{jk}$ is no larger than $O(n^{1/2})$. Alternatively, if $\omega_{jk} \to 0$, then (3.14) is minimized provided $m/n \to 0$. This discussion motivates an approach for choosing $m$ by tuning the bootstrapped probabilities to match the true probabilities. In reality, however, we do not know $\omega_{jk}$, $c$ or $Z$ so these must be estimated using $\hat{\omega}_{jk} = n^{1/2}(\hat{\theta}_k - \hat{\theta}_j)$, $\hat{c} = (\hat{\sigma}_j^2 + \hat{\sigma}_k^2)^{-1/2}$ and a random normal variable, respectively. The situation is simplified if we have a "gap" between the orders of the diverging $\omega_{jk}$ and those converging, such as the following:

(3.15)    For each pair $j, k$    either $\omega_{jk} \to 0$  or  $|\omega_{jk}|(\log n)^{-1/2} \to \infty$.

Thus, we estimate $m$ by choosing it to minimize the expression

(3.16)
$$\sum_{j,k;j \neq k} \int (\Phi[(m/n)^{1/2}\{-\hat{c}\hat{\omega}_{jk}(\log n)^{-1/2} + z\}]$$
$$- \Phi\{-\hat{c}\hat{\omega}_{jk}(\log n)^{-1/2}\})^2 \phi(z)\,dz.$$

The following theorem, a proof of which is given in the Ph.D. thesis of the second author, shows that choosing $m$ in this fashion is consistent.

THEOREM 3.2.  *Assume $p$ is fixed and that (3.1), (3.2), (3.5) and (3.15) hold. Choose $m$ by minimizing (3.16). Then we have for each $j$*

$$P(\hat{r}_j^* \leq r | \mathcal{X}) - P(\hat{r}_j \leq r) \to 0$$

*in probability.*

While this result suggests a way of determining $m$, there remains some uncertainty since the $(\log n)^{-1/2}$ factor used is not unique in generating good asymptotic performance. For example, replacing it with $(\log Cn)^{-1/2}$ for some constant $C$ would yield a similar theoretic result. In practice, the dataset under consideration often suggests whether the adopted factor is appropriate, and the choice of $m$ is reasonably robust against such changes.

3.5. *Theoretical properties in the case of large $p$.*  The results above can be generalized to cases where $p$ diverges with $n$ but the support of the limiting distribution of $\hat{r}_j$ remains bounded. The defining features of those extensions are that values of $|\theta_k - \theta_j|$, for indices $k$ that are not in the $\mathcal{K}_j$ of the previous section, should be at least as large as $(n^{-1}\log n)^{1/2}$; and values of $|\theta_k - \theta_j|$, for $k$ in $\mathcal{K}_j$, should be at least as small as $n^{-1/2}$. We shall give results of this type in Section 4.2. In the present section, we show how to capture, in a theoretical model, instances where both $p$ and the support of the distribution of $\hat{r}_j$ are large. Real-data examples of this type are given by Goldstein and Spiegelhalter (1996).

Specifically, we assume the following linear model for $\theta_j$:

(3.17)     $\theta_j = a - \varepsilon j$ for $1 \leq j \leq p$, where $a = a(n,p)$ does not depend on $j$ and $\varepsilon = \varepsilon(n) > 0$.

This condition ensures the simple numerical ordering $\theta_1 > \cdots > \theta_p$, which in more general contexts we can impose without loss of generality. Assumption (3.17) also allows us to adjust the difficulty of the empirical ranking problem by altering the size of $\varepsilon$; the difficulty increases as $\varepsilon$ decreases.

As in Theorem 3.1, we assume that the datasets $\mathcal{X}_j$ are independent, but now we permit $p = p(n)$ to diverge with $n$. In order that Theorem 3.3 below may be stated relatively simply, we assume that the quantities $Z_k = n^{1/2}(\hat{\theta}_k - \theta_k)$ all have the same asymptotic variance $\sigma$. Our main conclusion, that the standard $n$-out-of-$n$ bootstrap correctly captures order of magnitude but not constant multipliers, remains valid as long as the limiting variances of the $Z_k$'s are bounded away from zero and infinity.

We also assume conditions (3.18) and (3.19) below. In cases where each $\theta_j$ is a quantity such as a mean, a quantile, or any one of many different robust measures of location, those conditions follow from moderate-deviation properties of sums of independent random variables, provided the data have sufficiently many finite moments and $p$ does not diverge too rapidly as a function of $n$:

(3.18)     $P\{n^{1/2}(\hat{\theta}_k - \theta_k) \leq \sigma x\} = \Phi(x)\{1 + o(1)\} + o(p^{-1}n^{-1/2}\varepsilon^{-1})$, uniformly in $|x| = O(pn^{1/2}\varepsilon)$ and in $1 \leq k \leq p$, as $n \to \infty$, where $\sigma > 0$;

(3.19)     $P\{n^{1/2}(\hat{\theta}_k^* - \hat{\theta}_k) \leq \sigma x | \mathcal{X}\} = \Phi(x)\{1 + o_p(1)\} + o_p(p^{-1}n^{-1/2}\varepsilon^{-1})$, uniformly in $|x| = O(pn^{1/2}\varepsilon)$ and in $1 \leq k \leq p$, as $n \to \infty$, where $\sigma$ is as in (3.18).

In order for (3.18) and (3.19) to hold as $p$ increases, the value of $\varepsilon$ should decrease as a function of $p$, that is, the empirical ranking problem should be made more difficult for larger values of $p$. Define $\delta = (n/2)^{1/2}\varepsilon/\sigma$, where $\sigma > 0$ is as in (3.18) and (3.19), and put $\widetilde{\omega}_{jk} = n^{1/2}\{\hat{\theta}_k - \theta_k - (\hat{\theta}_j - \theta_j)\}/(2^{1/2}\sigma)$.

THEOREM 3.3.   *Assume that the datasets $\mathcal{X}_1, \ldots, \mathcal{X}_p$ are independent, that $\hat{\theta}_1^*, \ldots, \hat{\theta}_p^*$ are independent conditional on $\mathcal{X}$, that (3.17)–(3.19) hold, and that $p = p(n) \to \infty$ and $\varepsilon = \varepsilon(n) \downarrow 0$ as $n$ increases, in such a manner that $n^{1/2}\varepsilon \downarrow 0$ and $pn^{1/2}\varepsilon \to \infty$. Then*

(3.20)     $E(\hat{r}_j) = \delta^{-1} \int_{-j\delta}^{\infty} \Phi(-x)\,dx + o(\delta^{-1})$,

(3.21)     $E(\hat{r}_j^* | \mathcal{X}) = \{1 + o_p(1)\} \sum_{k\,:\,k \neq j} \Phi\{\widetilde{\omega}_{jk} + \delta(j - k)\} + o_p(\delta^{-1})$,

*uniformly in $1 \leq j \leq C/(n^{1/2}\varepsilon)$ for any $C > 0$.*

The implications of Theorem 3.3 can be seen most simply when $j$ is fixed, although other cases are similar. For any fixed $j$, it follows from (3.20) and (3.21) that

$$(3.22) \qquad E(\hat{r}_j) \sim C\delta^{-1},$$

$$(3.23) \qquad E(\hat{r}_j^*|\mathcal{X}) \sim_p \delta^{-1} \int_{-\infty}^{\infty} d\Phi(z) \int_0^{\infty} \Phi(W_j - x - z2^{-1/2}) \, dx,$$

where $C = \int_{x>0} \Phi(-x) \, dx$, $W_j = -(n/2)^{1/2}(\hat{\theta}_j - \theta_j)/\sigma$, $a_n \sim b_n$ for constants $a_n$ and $b_n$ means that $a_n/b_n \to 1$, and $A_n \sim_p B_n$ for random variables $A_n$ and $B_n$ means that $A_n/B_n \to 1$ in probability. Results (3.22) and (3.23) reflect the highly variable character of $\hat{r}_j$ in the difficult cases represented by the model (3.17). For example, if $r_j = j$, which of course is fixed if $j$ is fixed, then both $E(\hat{r}_j)$ and $E(\hat{r}_j^*|\mathcal{X})$ are of size $\delta^{-1}$, which diverges to infinity as $n \to \infty$. That is, despite $r_j$ being fixed, $\hat{r}_j$ tend to be so large that its expected value diverges. Similar arguments show that $\mathrm{var}(\hat{r}_j|\mathcal{X}_j)$ and $\mathrm{var}(\hat{r}_j^*|\mathcal{X}, \mathcal{X}_j^*)$ are both of size $\delta^{-1}$.

It is clear from (3.22) and (3.23) that the standard $n$-out-of-$n$ bootstrap correctly captures the order of magnitude, $\delta^{-1}$, of $E(\hat{r}_j)$, but does not get the constant multiplier right. Similar arguments, based on elementary properties of sums of independent random variables, show that the standard bootstrap produces a prediction interval for $\hat{r}_j$ for which the length has the correct order of magnitude, but again the constant multiplier is not correct. The $m$-out-of-$n$ bootstrap gets both the order of magnitude and the constant right, but at the expense of more restrictive conditions on $\varepsilon$; one could predict from Theorem 3.1 that this would be the case. It is also possible to establish a central limit theorem describing properties of $E(\hat{r}_j)$ and $E(\hat{r}_j^*|\mathcal{X})$. However, since the limitations of the bootstrap are clear at a coarser level than that type of analysis would address, then we shall not give those results here.

3.6. *Numerical properties.* We present numerical work which reinforces and complements the theoretical issues discussed above. In our first set of simulations, we observe $n$ independent data vectors $(X_1, \ldots, X_{10})$, where the $X_j$'s are independent and are, respectively, distributed as normal $N(\theta_j, 1)$. First, we consider the case where $\theta_j = 1 - (j/10)$, implying that the means are evenly spaced and do not depend on $n$. Although this model appears straightforward, the gaps between means are one tenth of the value of noise standard deviation, and so significant ranking challenges are present. However, Figure 1 shows that this is a case that the standard $n$-out-of-$n$ bootstrap can handle satisfactorily, with the 90% prediction intervals for the estimated ranks shrinking as $n$ grows.

Nevertheless, our theory suggests that the $n$-out-of-$n$ bootstrap will fail to correctly estimate the distribution in cases where the values of $\theta_j$ are relatively close. To investigate this issue, we took $\theta_j = 1$ for $j \in \{1, 2, 3, 4, 5\}$, and $\theta_j = 0$, otherwise. Then in our bootstrap replicates, we would expect $\hat{r}_j^*$, conditional on the data, to be approximately uniformly distributed on either the top five positions (in the case $j \leq 5$) or the bottom five (when $j \geq 6$). Figure 2 shows the difference in distributions for a simulation with $n = 1000$ and two choices of $m$. For each variable, the shading intensities in that column show the relative empirical distributions across ranks. Here the $m$-out-of-$n$ bootstrap, with $m = 300$, produces distributions closer to the truth, where each of the top-left and bottom-right regions would have exactly equal intensities everywhere.

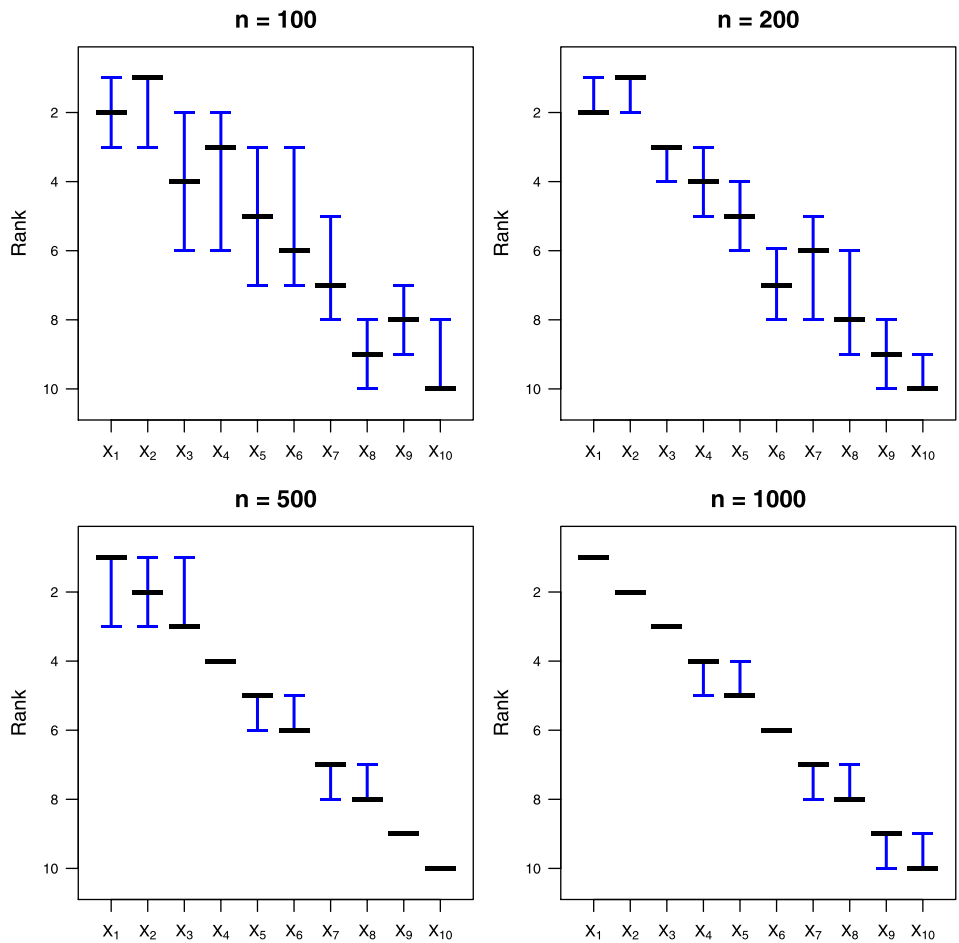

FIG. 1. *Ranking 90% prediction intervals for the case of fixed $\theta_j$.*

**n=1000, m=300**                    **n=1000, m=1000**

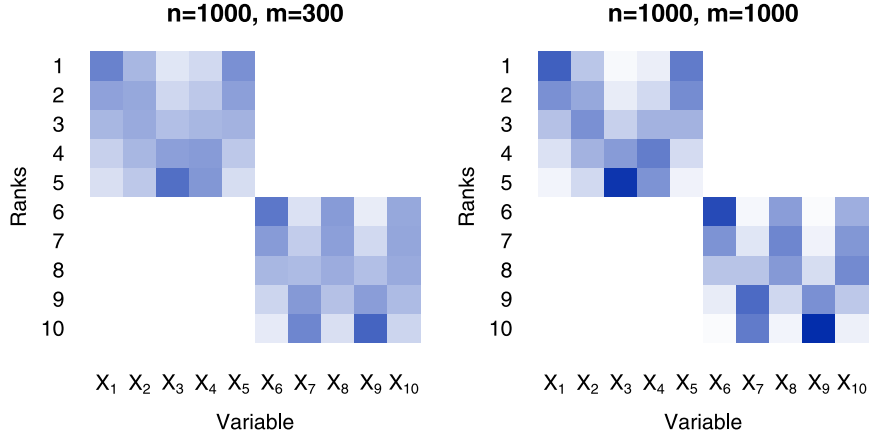

Fig. 2. *Distribution of ranks in the presence of ties.*

The case of perfect ties demonstrates the advantages of the $m$-out-of-$n$ bootstrap. In more subtle settings, when the $\theta_j$'s vary with $n$ and are not exactly tied, we are interested in the ability of the bootstrap to distinguish $\theta_j$'s for which the absolute differences in $|\theta_j - \theta_k|$ are relatively small. The theory suggests considering differences of size $m^{-\alpha}$, where $\alpha = \frac{1}{2}$ is the critical value, lower values of $\alpha$ tend toward a (degenerate) perfect separation of ranks, and higher values asymptotically behave as though $\theta_j$ and $\theta_k$ were tied. Therefore, the next set of simulations had the $\theta_j$'s equally spaced and uniformly decreasing, with $\theta_j - \theta_{j+1}$ equal to $0.2(10/m)^\alpha$. Here $m$ was taken to be $\min(10n^{1/2}, n)$. Figure 3 shows, for a given pair $(\alpha, n)$, the average number of ranks contained within the 90% rank prediction interval. The results accord with the theory; cases where $\alpha < 0.5$ tend toward perfect separation (an average of 1), and cases where $\alpha > 0.5$ tend toward completely random ordering (an average of 10). Situations where $\alpha$ is closer to 0.5 diverge more slowly, and the behavior when $\alpha = 0.5$ depends on the exact situation; in our simulations the degree of tuning has ensured that the case where $\alpha = 0.5$ does not show much tendency toward either extreme.

It is important to understand the distributional bias seen in the $n$-out-of-$n$ bootstrap. One way this can be done is by exploring the distribution implied by (3.7). The distribution is dependent on the realization of normal standard random variables $Z_1, \ldots, Z_p$. Figure 4 shows how the distribution of rankings varies with $Z_1$ for the special case of five variables, with $\omega_1 = \cdots = \omega_5 = 0$ and $Z_2 = Z_3 = Z_4 = Z_5 = 0$. Here, as $|Z_1|$ departs from 0, the ranking distribution is upset in two key ways. First, the average ranking is biased; for example, when $Z_1 = 1$ the average observed rank is 1.95 instead of 3, the average observed rank in the true underlying distribution obtained when $Z_1 = 0$. Second, the variation of the observed rank is reduced; the variance is 1.4 when $|Z_1| = 1$ compared with 2 in the true distribution. These two effects combine to give overconfidence in the $n$-out-of-$n$ bootstrap when it is not warranted.

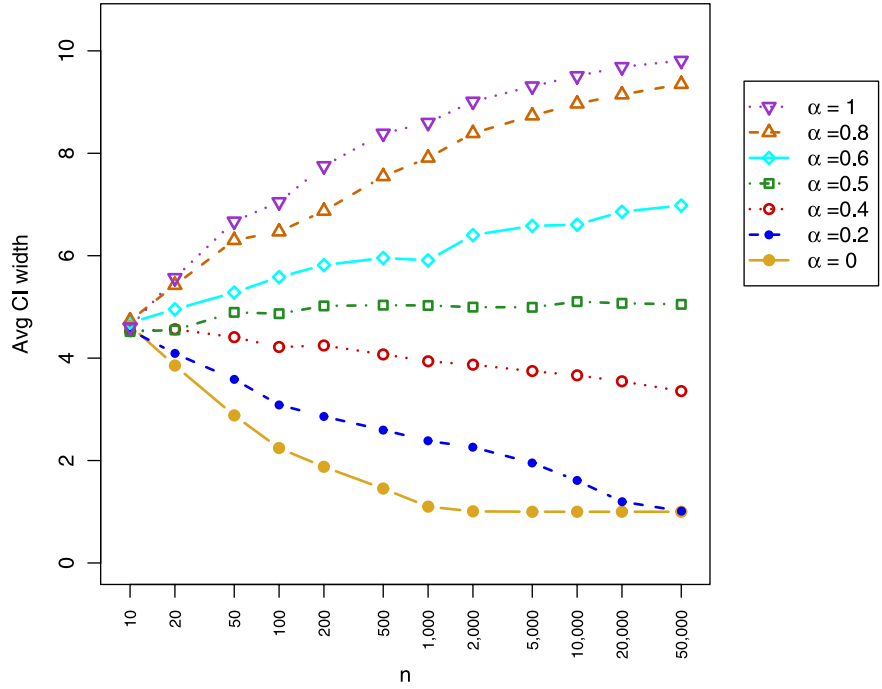

FIG. 3.   *Behavior of prediction interval widths for various* α.

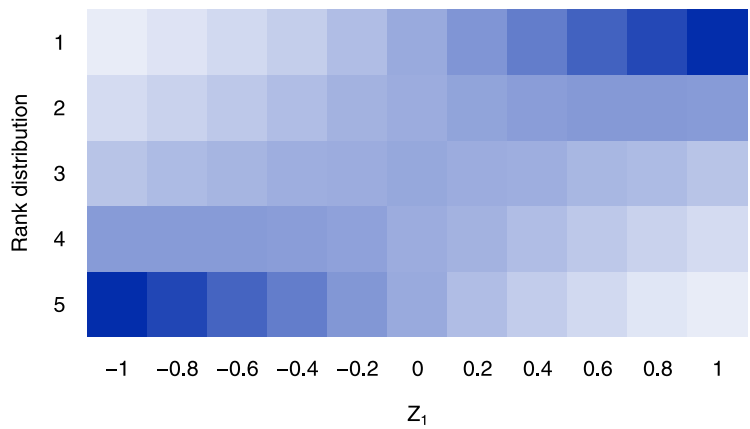

FIG. 4.   *Distribution of ranks for various* $Z_1$.

We now move to a real-data example. A service seeking to assist parents to choose secondary schools in the state of NSW, Australia, ranks 75 schools using the number of credits achieved in final year Higher School Certificate exams as a percentage of the number of exams sat. While there are clearly significant problems with such a simple statistic [see Goldstein and Spiegelhalter (1996)], the main one being that it ignores prior student ability, it would still be useful to give some indication of the variability of the rankings. Here $n_j$ represents the number of exams sat at school $j$, and $X_{ij}$, for $1 \leq i \leq n_j$, is an indicator variable for whether a credit was achieved in exam $i$. Then $\theta_j = E(X_{ij})$ and $\hat{\theta}_j = n_j^{-1} \sum_i X_{ij}$. Figure 5 shows 95% prediction intervals for the ranks using the $n$-out-of-$n$ bootstrap. It is clear that caution needs to be exercised when interpreting the intervals, the average width of which exceeds 14 places. However, we know that the $n$-out-of-$n$ bootstrap ranking understates the true uncertainty, which would be better captured using the $m$-out-of-$n$ bootstrap. Figure 6 shows the results using $m_j = \lfloor n_j \times 35.5\% \rfloor$. The percentage here was chosen using the approach discussed in Section 3.4, attempting to minimize the squared error between the bootstrap and real ranking distributions. Observe that the widths of the prediction intervals are now markedly longer (58% longer on average); the widest confidence interval now covers 81% of the possible rankings. Our theoretical results argue that these longer widths give a better indication of the true uncertainty associated with the ranking.

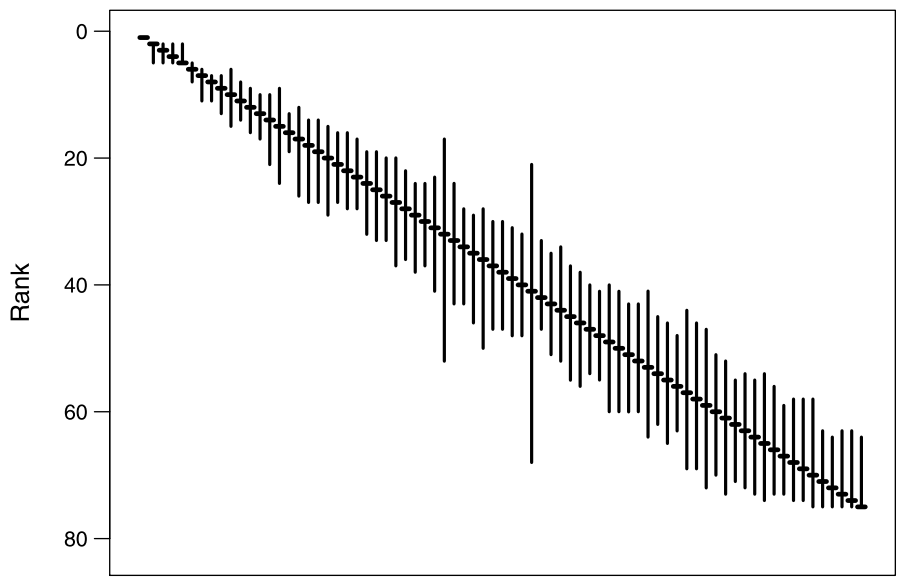

Fig. 5. *School ranking prediction intervals for n-out-of-n bootstrap.*

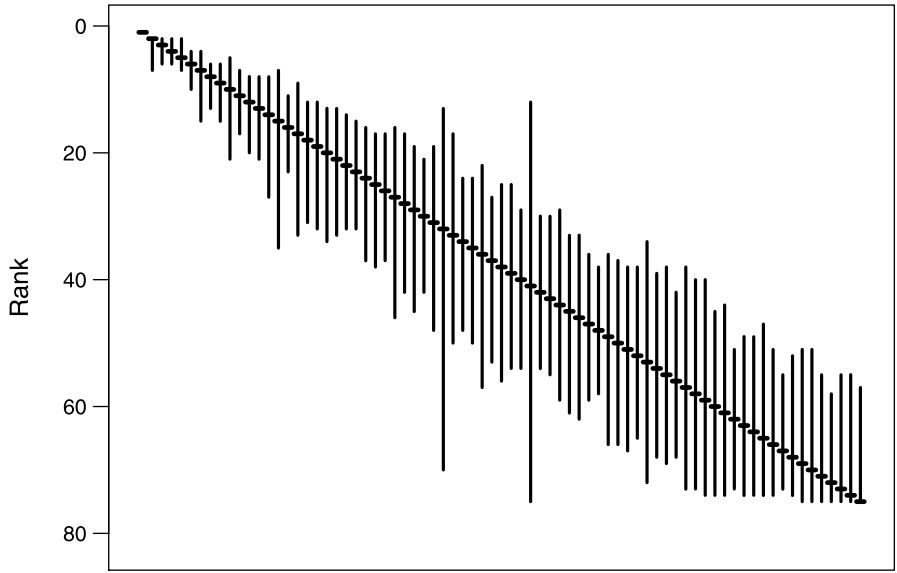

Fig. 6. *School ranking prediction intervals for m-out-of-n bootstrap with $m_j$ equal to 35.5% of $n_j$.*

## 4. Properties in cases where the data come as independent $p$-vectors.

4.1. *Motivation for the independent-component bootstrap.* In this section, we argue that when vector components are not strongly dependent the standard, "synchronous" bootstrap may distort relationships among components, particularly in the setting of conditional inference and when $p$ is large. In such cases, even if the assumption of independent components is not strictly correct, it may be advantageous to apply the bootstrap as though independence prevailed. We refer to this working assumption as that of "component-wise independence."

We treat the case where the data arise via a sample $\mathcal{X}_0 = \{X_1, \ldots, X_n\}$ of independent $p$ vectors. Here $X_i = (X_{i1}, \ldots, X_{ip})$, and $\mathcal{X}_j = \{X_{1j}, \ldots, X_{nj}\}$ denotes the set of $j$th components. The conventional, synchronous form of the nonparametric bootstrap involves the following resampling algorithm:

(4.1)   Draw a resample $\mathcal{X}_0^* = \{X_1^*, \ldots, X_m^*\}$ by sampling randomly, with replacement, from $\mathcal{X}_0$, write $X_i^* = (X_{i1}^*, \ldots, X_{ip}^*)$ and take $\mathcal{X}_j^* = \{X_{1j}^*, \ldots, X_{mj}^*\}$.

We can view $\mathcal{X}_j^*$ as the resample drawn from the $j$th "population." In (4.1), we take $m \leq n$, thereby allowing for the $m$-out-of-$n$ bootstrap.

We argue that this bootstrap method is not always satisfactory in problems where ranking is involved. One reason is that:

(4.2)   If the data have a continuous distribution, then knowing the dataset $\mathcal{X}_j^*$ conveys perfect information about which data vectors $X_i$ are included in $\mathcal{X}_0^*$, defined in (4.1), and with what frequencies. Hence, knowing $\mathcal{X}_j^*$ tells us $\mathcal{X}_k^*$ for each $k$, and in particular the resamples $\mathcal{X}_1^*, \ldots, \mathcal{X}_p^*$ cannot be regarded as independent, conditional on $\mathcal{X}_0$, even if the vector components are independent.

This result holds for the $m$-out-of-$n$ bootstrap as well as for the standard, synchronous bootstrap, and so the problems to which it leads cannot be alleviated simply by passing to a smaller resample size.

To elucidate the consequences of (4.2), note that the $j$th empirical rank $\hat{r}_j$, and its bootstrap version $\hat{r}_j^*$, can be written as

$$(4.3) \qquad \hat{r}_j = 1 + \sum_{k:\, k \neq j} I(\hat{\theta}_j \leq \hat{\theta}_k), \qquad \hat{r}_j^* = 1 + \sum_{k:\, k \neq j} I(\hat{\theta}_j^* \leq \hat{\theta}_k^*),$$

respectively. Here, $\hat{r}_j$ and $\hat{r}_j^*$ are as at (2.2) and (2.3). We wish to estimate aspects of the distribution of $\hat{r}_j$. For example, we might seek an estimator of the variance of the conditional mean, $\hat{u}_j = E(\hat{r}_j | \mathcal{X}_j) = 1 + \sum_{k:\, k \neq j} \pi_{jk}$, of $\hat{r}_j$ given $\mathcal{X}_j$; or we might wish to approximate the variance of $\hat{r}_j$. [To derive the formula for $\hat{u}_j$, we used the first part of (4.3), and took $\pi_{jk} = P(\hat{\theta}_j \leq \hat{\theta}_k | \mathcal{X}_j)$.] Undertaking conditional inference is particularly attractive in problems where $p$ is large, because it has the potential to greatly reduce variability, from $O(p^2)$ (the order of the unconditional variance of $\hat{r}_j$) to $O(p)$ (the order of the variance of $\hat{r}_j$, conditional on $\mathcal{X}_j$, if the components are sufficiently weakly dependent).

The bootstrap version of $\hat{u}_j$ can be computed using the second formula in (4.3): $\hat{u}_j^* = E(\hat{r}_j^* | \mathcal{X}, \mathcal{X}_j^*) = 1 + \sum_{k:\, k \neq j} \hat{\pi}_{jk}^*$, where $\mathcal{X} = \bigcup_k \mathcal{X}_k$ and $\hat{\pi}_{jk}^* = P(\hat{\theta}_j^* \leq \hat{\theta}_k^* | \mathcal{X}, \mathcal{X}_j^*)$. If we use the synchronous bootstrap algorithm at (4.1), then it follows from (4.2) that $\hat{\pi}_{jk}^* = I(\hat{\theta}_j^* \leq \hat{\theta}_k^*)$. Since the probability has degenerated to an indicator function, then even when using the $m$-out-of-$n$ bootstrap, and in the conventional setting of fixed $p$ and increasing $n$, $\mathrm{var}(\hat{u}_j^* | \mathcal{X}) - \mathrm{var}(\hat{u}_j)$ fails to converge to zero except in degenerate cases.

The errors can become still more pronounced if $p$ diverges with $n$. Indeed, in the problem of estimating

$$\mathrm{var}(\hat{u}_j) = \sum_{k_1:\, k_1 \neq j} \sum_{k_2:\, k_2 \neq j} \mathrm{cov}(\pi_{jk_1}, \pi_{jk_2})$$

using

$$(4.4) \qquad \mathrm{var}(\hat{u}_j^* | \mathcal{X}) = \sum_{k_1:\, k_1 \neq j} \sum_{k_2:\, k_2 \neq j} \mathrm{cov}(\pi_{jk_1}^*, \pi_{jk_2}^* | \mathcal{X})$$

and in the context of component-wise independence, the synchronous bootstrap at (4.1) introduces correlation terms of size $n^{-1/2}, n^{-1}, \ldots$; those terms would be zero if the bootstrap algorithm correctly reflected component-wise independence. If $p$ is much larger than $n$, then the impact of the extraneous terms is magnified by the summations over $k_1$ and $k_2$ in (4.4). These problems, too, persist when employing the $m$-out-of-$n$ bootstrap.

The situation improves significantly if, instead of using the synchronous bootstrap at (4.1), we employ the following independent-component resampling algorithm:

(4.5)    Compute $\mathcal{X}_j^* = \{X_{1j}^*, \ldots, X_{mj}^*\}$ by sampling randomly, with replacement, from $\mathcal{X}_j = \{X_{1j}, \ldots, X_{nj}\}$; and do this independently for each $j$.

In this case, when using the $m$-out-of-$n$ bootstrap and working under the assumption of component-wise independence, $\mathrm{var}(\hat{u}_j^* | \mathcal{X}) - \mathrm{var}(\hat{u}_j)$ converges to zero as $n$ diverges, and the bothersome $n^{-1/2}$ terms that arise when estimating $\mathrm{var}(\hat{u}_j)$, using the synchronous bootstrap, vanish. To summarize, under component-wise independence the independent-component bootstrap, defined at (4.5), corrects for significant errors that can be committed by the synchronous bootstrap algorithm at (4.1).

Importantly, similar conclusions are also reached in cases where $p$ is large and the component vectors $(X_{1j}, \ldots, X_{nj})$ are not independent. In particular, if the dependence among components is sufficiently weak to ensure that the asymptotic distribution of $\hat{r}_j$ is identical to what it would be if the components were independent, then the independent-component bootstrap has obvious attractions. For example, in inferential problems involving conditioning on $\mathcal{X}_j$, it gives statistical consistency in contexts where the synchronous bootstrap does not. This can happen even under conditions of reasonably strong dependence, simply because the highly ranked components are lagged well apart. Details will be outlined in the first paragraph of Section 4.3.

4.2. *Theoretical properties.* We address only the $j_0$ highest-ranked populations, which for notational convenience we take to be those with indices $j = 1, \ldots, j_0$, and we take the ranks of these populations to be virtually tied, so that the limiting distribution of $\hat{r}_j$ is nondegenerate. Also, we allow both $p$ and the distributions of $\Pi_1, \ldots, \Pi_p$ to depend on $n$. In particular, we assume that:

(4.6)       $n^{1/2}(\theta_1 - \theta_j) \to 0$     for $j = 1, \ldots, j_0$,

(4.7)              $p = o(n^{C_1})$     for some $C_1 > 0$.

To determine the limiting distribution of $\hat{r}_j$, we further suppose that

$$(4.8) \qquad (n/\log n)^{1/2} \inf_{j_0 < j \le p}(\theta_1 - \theta_j) \to \infty$$

and

(4.9)  the random variables $n^{1/2}(\hat{\theta}_j - \theta_j)$, for $1 \le j \le j_0$, are asymptotically independent and normally distributed with zero means and respective variances $\sigma_j^2$; and, for $C_2 > 0$ sufficiently large, $\sup_{j \le p} P\{|\hat{\theta}_j - \theta_j| > C_2(n^{-1} \log n)^{1/2}\} = O(n^{-C_1})$.

When discussing the efficacy of the $m$-out-of-$n$ bootstrap we ask, instead of (4.6), (4.8) and (4.9), that

$$(4.10) \qquad m^{1/2}(\theta_1 - \theta_j) \to 0 \qquad \text{for } j = 1, \ldots, j_0,$$

$$(4.11) \qquad (m/\log m)^{1/2} \inf_{j_0 < j \le p}(\theta_1 - \theta_j) \to \infty,$$

(4.12)  conditional on $\mathcal{X}$, the random variables $m^{1/2}(\hat{\theta}_j^* - \hat{\theta}_j)$, for $1 \le j \le j_0$, are asymptotically independent and normally distributed with zero means and respective variances $\sigma_j^2$; and for $C_2 > 0$ sufficiently large, $\sup_{j \le p} P\{|\hat{\theta}_j^* - \hat{\theta}_j| > C_2(m^{-1} \log m)^{1/2}\} = O(n^{-C_1})$.

For example, the last parts of (4.9) and (4.12) hold if $\theta_j$ and $\hat{\theta}_j$ are, respectively, population and sample means, if the associated population variances are bounded away from zero, and if the supremum over $j$ of absolute moments of order $C_3$, for the population $\Pi_j$, is bounded for a sufficiently large $C_3 > 0$. See Rubin and Sethuraman (1965) and Amosova (1972). Likewise, (4.9) and (4.12) also apply in cases where each $\theta_j$ is a quantile or any one of many different robust measures of location. The first part of (4.9) is a standard central limit theorem for the estimators $\hat{\theta}_j$, and so is a weak assumption. In (4.12), we do not specify using the independent-component bootstrap [see (4.5)], but if we do impose that condition then the first part of (4.12) is a conventional central limit theorem for the $m$-out-of-$n$ bootstrap, and in that setting we do not need to assume independence of the asymptotic normal distribution of the variables $m^{1/2}(\hat{\theta}_j^* - \hat{\theta}_j)$; it follows from the nature of the independent-component bootstrap.

THEOREM 4.1. *Let* $1 \le j \le j_0$.

(i) *If* (4.6)–(4.9) *hold then the ranks* $\hat{r}_1, \ldots, \hat{r}_{j_0}$ *are asymptotically jointly distributed as* $R_1, \ldots, R_{j_0}$, *where*

$$(4.13) \qquad R_j = 1 + \sum_{k \,:\, k \le j_0, k \ne j} I(Z_j \sigma_j \le Z_k \sigma_k)$$

*and* $Z_1, \ldots, Z_{j_0}$ *are independent and normal* $N(0, 1)$.

(ii) *Assume (4.7) and (4.9)–(4.12), and use the m-out-of-n bootstrap (where $m/n \to 0$ and $m \to \infty$ as $n \to \infty$), in either the conventional form at (4.2) or the component-wise from at (4.5). Then the distribution of $(\hat{r}_1^*, \ldots, \hat{r}_{j_0}^*)$, conditional on the data $\mathcal{X}$, converges in probability to the distribution of $(R_1, \ldots, R_{j_0})$.*

(iii) *Assume (4.7) and (4.9)–(4.12), use the m-out-of-n bootstrap with $m/n \to 0$ and $m \to \infty$ as $n \to \infty$, and implement the bootstrap component-wise, as in (4.5). Then the distribution of $\hat{u}_j^*$, conditional on $\mathcal{X}$, is consistent for that of $\hat{u}_j$. That is,*

$$(4.14) \qquad P\{E(\hat{r}_j^* | \mathcal{X}, \mathcal{X}_j^*) \le x | \mathcal{X}\} \to P\{E(R_j | Z_j) \le x\}$$

*in probability, for all continuity points $x$ of the cumulative distribution function $P\{E(R_j | Z_j) \le x\}$. Moreover, $\mathrm{var}(\hat{r}_j^* | \mathcal{X}) \to \mathrm{var}(R_j)$.*

4.3. *Discussion.* The assumptions underpinning Theorem 4.1 do not require the components of the data vectors $X_i = (X_{i1}, \ldots, X_{ip})$ to be independent, but they do ask that the empirical ranks $\hat{\theta}_j$, corresponding to the true $\theta_j$'s that are virtually tied for the top $j_0$ positions, be asymptotically independent. See the first part of (4.9). That condition holds in many problems where $p$ is diverging but the components are strongly dependent, for example, when $\theta_j$ is a mean and the common distribution of the vectors $X_i$ is determined by adding $\theta_j$'s randomly to centred, although potentially strongly dependent, noise. For example, if the components of the noise process are $q$-dependent, where the integer $q$ is permitted to diverge with increasing $n$ and $p$, then in the case of fixed $j_0$ explored in Theorem 4.1, sufficient independence is ensured by the condition that $q/p \to 0$ as $p \to \infty$.

Parts (i) and (ii) of Theorem 4.1 together imply that (3.13) continues to hold in the present setting, provided $j$ is in the range $1 \le j \le j_0$.

As noted in Section 4.1, the result in the first part of Theorem 4.1(iii) does not hold if the synchronous bootstrap is used. Likewise, while the independent-component, $m$-out-of-$n$ bootstrap can be proved to consistently estimate the distribution of $\mathrm{var}(\hat{r}_j | \mathcal{X}_j)$, neither the $n$-out-of-$n$ bootstrap nor its $m$-out-of-$n$ bootstrap form give consistency if applied using the conventional resampling algorithm at (4.1). The same challenges arise for a variety of other estimation problems; the problems treated in Theorem 4.1(iii) are merely examples.

In cases where $p$ is very much larger than $n$, and the aim is to discover information concealed in a very high-dimensional dataset, choosing $m$ for the $m$-out-of-$n$ bootstrap might best be regarded as selecting the level of sensitivity rather than as choosing the level of smoothing in a more conventional, $m$-out-of-$n$ bootstrap sense. Since the desired level of sensitivity depends on the unknown populations $\Pi_j$, and, in the most important marginal cases, is unknown, then it may not always be appropriate to use a standard empirical approach to choosing $m$. Instead, numerical results for different values of $m$ could be obtained.

Results analogous to Theorem 3.3 can also be established in the present setting. In particular, in cases where $\hat{r}_j$ is highly variable, the standard $n$-out-of-$n$ bootstrap correctly captures the order of magnitude, but not the constant multiplier, of characteristics of the distribution of $\hat{r}_j$, for example, its expected value and the lengths of associated prediction intervals.

4.4. *Numerical properties.*  To gain insight into the advantages of the independent-component bootstrap we consider the following setting. Suppose we have $p$ variables and $n$ observations, and the $j$th variable $X_j$ is modelled by $X_j = \theta_j + Z_j$, where $\theta_j$ is a constant, $Z_j$ is a standard random normal variable, $\mathrm{cor}(Z_j, Z_k) = \rho_n$ when $j \neq k$, and $\rho_n$ decreases to 0 as $n$ increases. We wish to compare performance of the standard and independent-component bootstraps in the task of ranking the values of $\theta_j$. As our performance measure, we use the squared error criterion:

$$\sum_j \sum_r E\{P(\hat{r}_j^* = r | \mathcal{X}) - P(\hat{r}_j = r)\}^2.$$

Figure 7 gives results for $n = 50$, $p = 200$ and $\theta_j = 1 - \{j/(p-1)\}$, for various choices of $\rho_n$. It shows that the independent-component bootstrap consistently improves performance. Interestingly, performance of the independent-component case is at its best when a reasonable level of correlation present. This is apparently because, in the presence of correlation, the true ranking distribution becomes more "lumpy" or more degenerate.

The Ro131 dataset was used by Segal et al. (2003) to compare a variety of genomic approaches. The dataset contains 30 independent observations, each with continuous expression levels $X_i$ for 6,319 genes, as well as a continuous response $Y_i$ for the expression level of a G protein-coupled receptor called Ro131. Hall and Miller (2009) used generalized correlation between the observed $Y$ and each set of gene expressions $\mathcal{X}_j$ to rank the genes, and then applied the standard, synchronous bootstrap to give indicative prediction intervals for these rankings. These results are presented in Figure 8 for the top 15 variables. It should be observed that significant levels of correlation exist between pairs of influential genes. There are at least two possible reasons for this. First, if gene expression levels closely follow the movements of response variables then genes will share some of this correlation indirectly. Second, there may be intrinsic correlation between two genes if they are controlled by some common underlying process.

If the first reason is suspected to be the dominant one, then the independent-component bootstrap should give a better indication of uncertainties in ranking. Figure 9 depicts results for the independent-component bootstrap.

Notice that prediction interval widths are greater than in the synchronous case. This is because the positive correlations among values of $\hat{\theta}_j$ in the synchronous case reduce the variations in rankings.

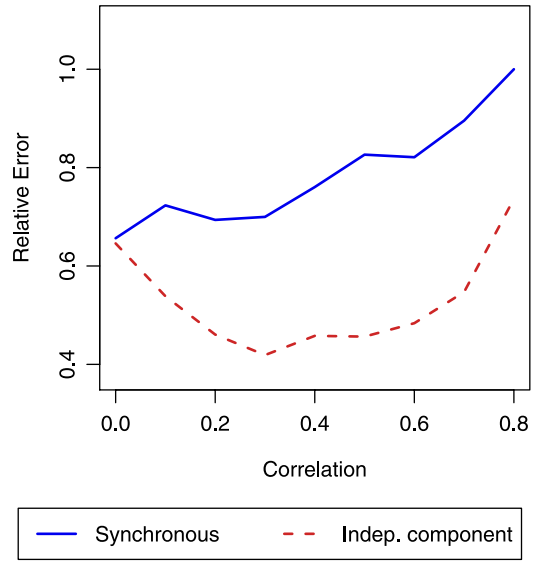

FIG. 7. *Relative error of synchronous and independent-component bootstrap distributions.*

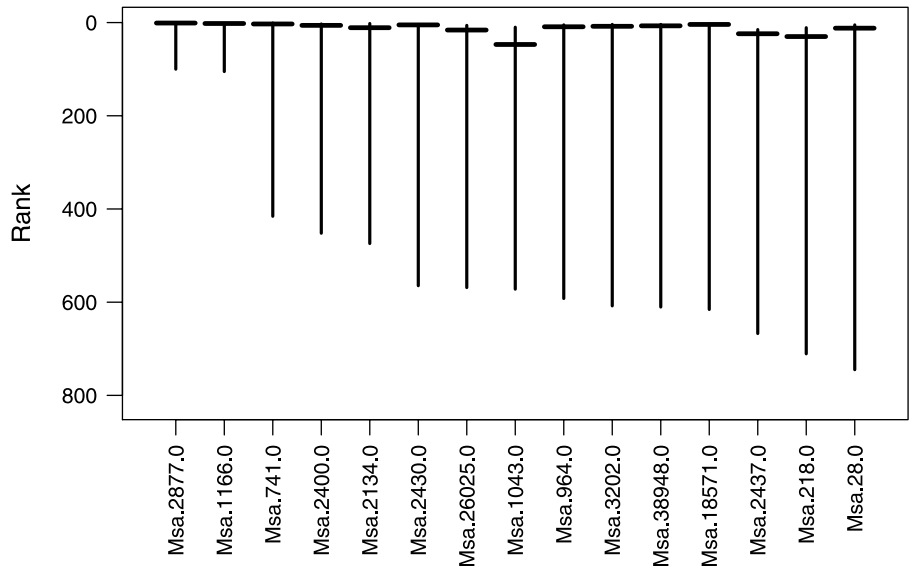

FIG. 8. *Synchronous bootstrap results for Ro131 dataset.*

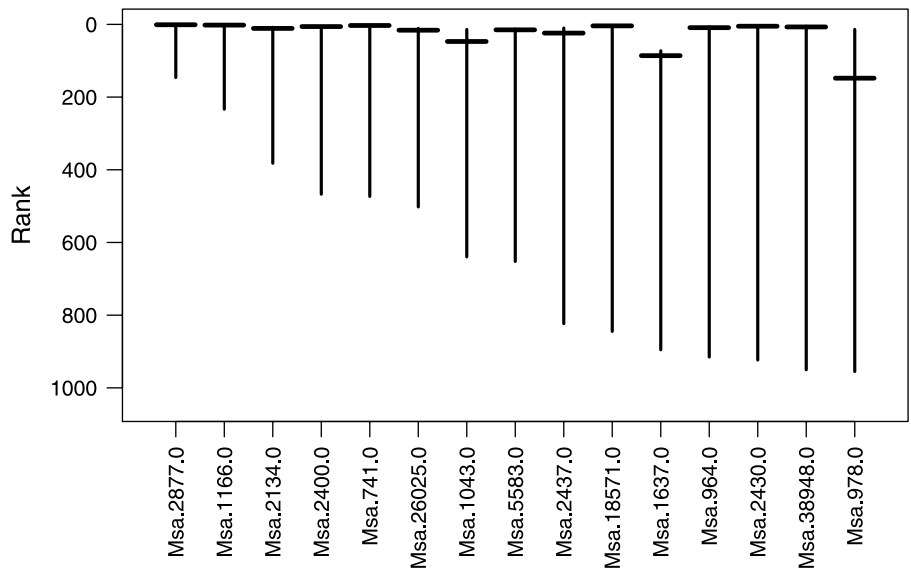

Fig. 9. *Independent-component bootstrap results for Ro131 dataset.*

Another plot that is useful in understanding rankings is that of conditional rankings, the subject of Theorem 4.1. Figure 10 shows the rankings for the top genes, together with prediction intervals for $\hat{r}_j^*$, conditional on both $\mathcal{X}$ and $\hat{\theta}_j^*$. Thus, for a given gene we have held the observed generalised correlation for it constant and bootstrapped on all other genes, to estimate how the genes should be ranked given the value of $\hat{\theta}_j$. The results for this analysis are highly dependent on whether the bootstrap is performed synchronously or independently. For reasons given in Sections 4.1–4.3, we prefer the independent-component bootstrap in this situation. Figure 10 displays the corresponding prediction intervals. Two features of the results are striking. First, the prediction intervals are very narrow compared to those seen in Figures 8 and 9, highlighting that the fact that most of the uncertainty in ranking the $j$th gene comes from the uncertainty of $\hat{\theta}_j$ itself. Second, the prediction intervals lie below the actual point estimate for the rank. This suggests that if the experiment were performed again, we would be unlikely to see the top-ranked variables rank as highly as before. In fact, we would expect the top variable to rank outside the top twenty, even if it appeared as strongly as it did in our observed data. These two observations are interesting, and highlight the challenges of variable selection in such high-dimensional settings.

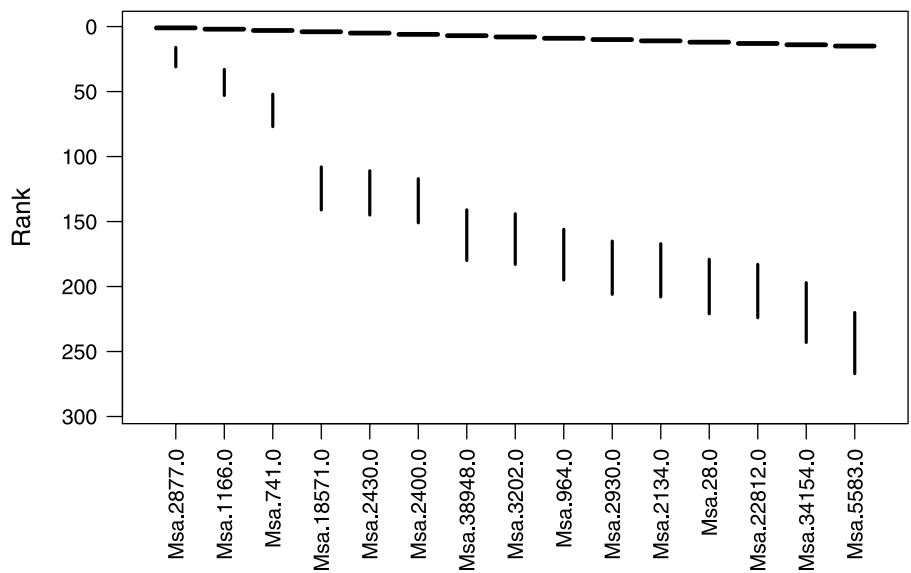

Fig. 10. *Independent reverse synchronous bootstrap results for Ro131 dataset.*

We reiterate here one observation relevant to both the independent-component bootstrap and the discussion of the $m$-out-of-$n$ bootstrap in Section 3. When constructing prediction intervals for ranks, the method that produces the shortest intervals is not necessarily the most powerful or the most accurate. Both the theoretical and numerical results suggest that the synchronous bootstrap will produce widths that are too narrow compared to the theoretical ranking distribution; the bootstrap ranks become "anchored" to the observed empirical ranks. Thus, interpreting ranking sensitivities for real datasets involves attempting to balance both maximising the power of an approach with the risks of overstating ranking accuracy. In many cases, it will be simulation and experimentation that suggest the best balance in a given situation.

The final example comprises of a set of simulations that illustrate the results of Theorem 4.1 in a high-dimensional setting. The aim here is to estimate the correct distribution for the top five ranked variables. For each of six scenarios, we start with the base case of $n = 20$, $p = 500$, which was constructed as follows. The mean is once more the statistic of interest. Each data point $X_{ij}$ is normal with standard deviation 0.25 and the $j$th mean is $\theta_j = 1$ for $j = 1, \ldots, 5$ and is randomly sampled from the uniform distribution over $[0, 0.9]$ when $j > 5$. Once the data is generated, we may derive the ranking distribution using the independent component bootstrap with $m = 20$. We use the statistic

$$\text{Error} = \sum_{j=1}^{5} \sum_{r=1}^{5} \{P(\hat{r}^* \leq r | \mathcal{X}) - P(\hat{r} \leq r)\}^2,$$

to measure how accurately the rankings for the first five variables are estimated. Notice that $P(\hat{r} \leq r) = r/5$ for $r = 1, \ldots, 5$ and that the error statistic is 0 if and only if this distribution is matched exactly in the bootstrapped distribution. We repeat this experiment 100 times and report the average error along with 90% confidence intervals for this average. From here, the simulation grows by increasing $n$ and increasing $m$ at rate $n/\log(n)$. In each scenario, $p$ is constant or grows at a linear or quadratic rate relative to $n$. Also, the gap between the mean of the top five variables and the upper range of the uniform sampling distribution is either left constant or shrunk at a square rooted logarithmic rate. This results in six scenarios, the results of which are plotted in Figure 11. The error has been scaled so that 100 denotes maximum possible error. Observe that the quadratic growth simulations in particular achieve very high dimensions; when $n = 140$, $p = 24,500$, which is competitive with the dimensionality for many genomic applications.

Theorem 4.1 establishes that under each of these scenarios the distribution of the top five variables should be estimated correctly, since $p$ increases only polynomially and the gap is either constant or shrinks sufficiently slowly; compare with (4.7), (4.12). The results reinforce these findings, with error steadily decreasing in all cases except the quadratic ones. In these final cases, the error increases briefly until the stability of the means outweighed the effects of increasing $p$ and decreasing gap. The error then steadily decreases, albeit at a much slower rate than the constant and linear scenarios. We can see that the problem is noticeably more difficult when the gap shrinks, as well as when $p$ grows at a faster rate.

This example was constructed to demonstrate that the theoretical results can hold while the data size remained computable. However, there are instances where very large $n$ are needed before such distributional accuracy is obtained. For instance, if we tripled the standard deviation in final scenario, where we have quadratic growth in $p$ and a shrinking gap, we would require $n > 1000$ before the error started to decrease and satisfactory results were obtained. In this case, $p$ would be over one million, which is in excess of current desktop computer capability.

## 5. Technical arguments.

5.1. *Proof of Theorem 3.1.* (i) In view of the first part of (3.5), we may write

$$(5.1) \qquad \hat{r}_j = 1 + \sum_{k:\, k \neq j} I(\hat{\theta}_j \leq \hat{\theta}_k) = 1 + \sum_{k:\, k \neq j} I(\sigma_j \Delta_j \leq \sigma_k \Delta_k + \omega_{jk}),$$

where the random variables $\Delta_k = n^{1/2}(\hat{\theta}_k - \theta_k)/\sigma_k$ are jointly independent and asymptotically standard normal. Result (3.6) can be proved from this quantity by considering the respective cases where values in the sequence $\omega_{jk}$, for $1 \le k \le p$, are finite or infinite.

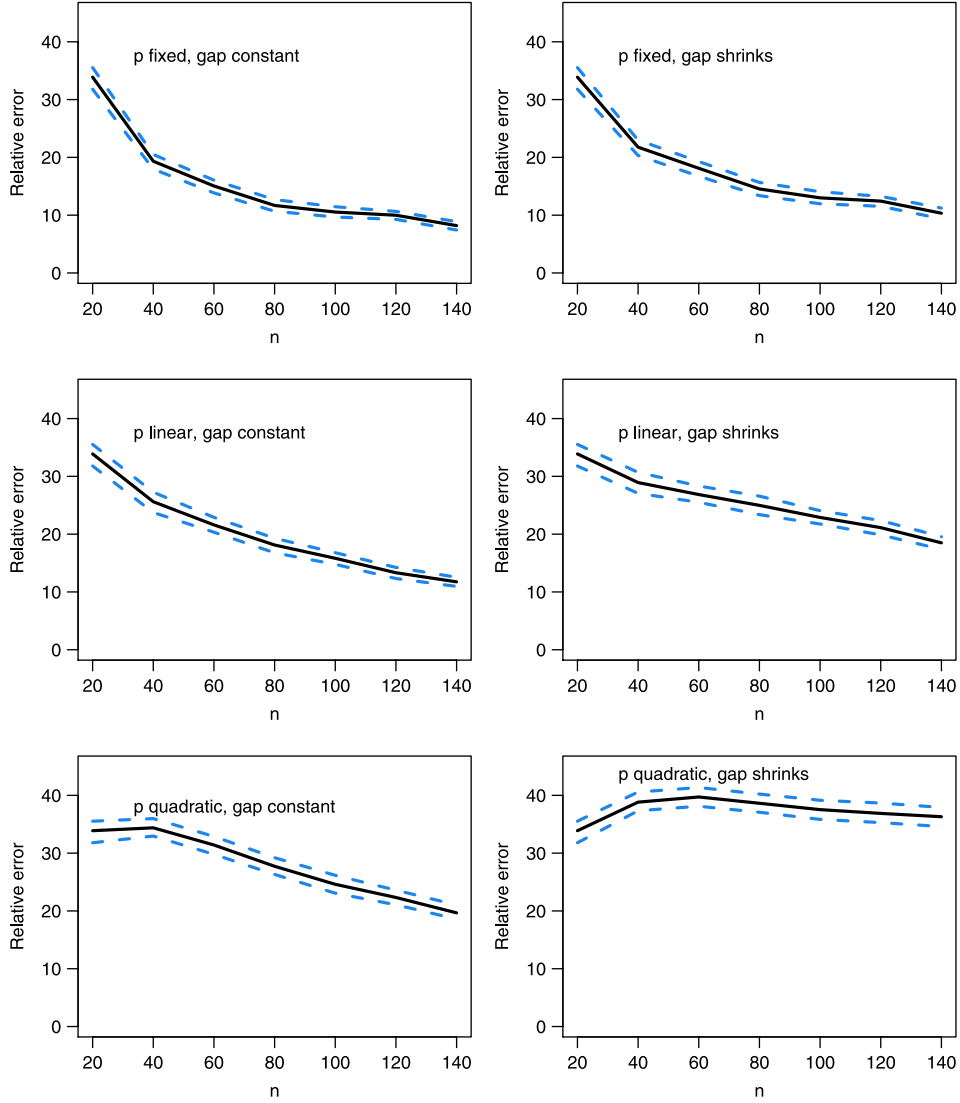

Fig. 11.   *Average error with 90% confidence intervals for $p > n$ simulations.*

(ii) To derive (3.7), we note that, in view of the second part of (3.5),

$$\hat{r}_j^* = 1 + \sum_{k\,:\,k \ne j} I(\hat{\theta}_j^* \le \hat{\theta}_k^*)$$

$$(5.2) \qquad = 1 + \sum_{k\,:\,k \ne j} I(n^{-1/2}\sigma_j\Delta_j^* \le n^{-1/2}\sigma_k\Delta_k^* + \hat{\theta}_k - \hat{\theta}_j)$$

$$= 1 + \sum_{k \,:\, k \neq j} I(\sigma_j \Delta_j + \sigma_j \Delta_j^* \leq \sigma_k \Delta_k + \sigma_k \Delta_k^* + \omega_{jk}),$$

where, conditional on $\mathcal{X}$, the random variables $\Delta_k^* = n^{1/2}(\hat{\theta}_k^* - \hat{\theta}_k)/\sigma_k$ are jointly independent and asymptotically standard normal, and the $\Delta_k$'s are as in (5.1). Since, by the first part of (3.5), the $\Delta_k$'s are asymptotically independent and standard normal (in an unconditional sense), then by Kolmogorov's extension theorem, we can (on a sufficiently large probability space) find random variables $Z_1, \ldots, Z_p$ which depend on $n$, are exactly independent and exactly standard normal for each $n$, and have the property that $\Delta_k = Z_k + o_p(1)$ for each $k$, as $n \to \infty$. Result (3.7) follows from these properties and (5.2).

(iii) Result (5.2) continues to hold in the case of the $m$-out-of-$n$ bootstrap, except that to obtain the arguments of the indicator functions there we have to multiply throughout by $m^{1/2}$ rather than $n^{1/2}$. This means that to interpret (5.2), we should redefine $\Delta_k = m^{1/2}(\hat{\theta}_k - \theta_k)/\sigma_k$ and $\Delta_k^* = m^{1/2}(\hat{\theta}_k^* - \hat{\theta}_k)/\sigma_k$. Since $m/n \to 0$, then, on the present occasion, $\Delta_k \to 0$ in probability for each $k$, but, in view of the second part of (3.5), the conditional distribution of $\Delta_k^*$ continues to be asymptotically normal $N(0, \sigma_k^2)$. Result (3.8) now follows from (5.2).

5.2. *Proof of Theorem 3.3.* Observe from (3.17), (3.18) and (5.1) that

$$E(\hat{r}_j) - 1 = \sum_{k \,:\, k \neq j} P(\hat{\theta}_j \leq \hat{\theta}_k) = \sum_{k \,:\, k \neq j} P\{\Delta_j \leq \Delta_k + 2^{1/2}\delta(j-k)\}$$

$$= \{1 + o(1)\} \sum_{k \,:\, k \neq j} \Phi\{\delta(j-k)\} + o(\delta^{-1})$$

$$= \delta^{-1} \int_{-j\delta}^{\infty} \Phi(-x)\,dx + o(\delta^{-1}),$$

where $\Delta_k = n^{1/2}(\hat{\theta}_k - \theta_k)/\sigma$. This gives (3.20). Similarly, (3.21) follows from

$$E(\hat{r}_j^* | \mathcal{X}) - 1 = \sum_{k \,:\, k \neq j} P(\hat{\theta}_j^* \leq \hat{\theta}_k^* | \mathcal{X})$$

$$= \sum_{k \,:\, k \neq j} P\{\Delta_j^* \leq \Delta_k^* + 2^{1/2}\tilde{\omega}_{jk} + 2^{1/2}\delta(j-k)\}$$

$$= \{1 + o_p(1)\} \sum_{k \,:\, k \neq j} \Phi\{\tilde{\omega}_{jk} + \delta(j-k)\} + o_p(\delta^{-1}),$$

where $\Delta_k^* = n^{1/2}(\hat{\theta}_k^* - \hat{\theta}_k)/\sigma$.

5.3. *Proof of Theorem 4.1.* (i) By (4.9), the probability that $|\hat{\theta}_j - \theta_j| > C_2(n^{-1}\log n)^{1/2}$ for some $j = 1, \ldots, p$, equals $O(pn^{-C_1}) = o(1)$, where we used (4.7) to obtain the last identity. Therefore, by (4.6), (4.9) and (4.8), for each $C > 0$, the probability that $\hat{\theta}_j - \hat{\theta}_k > C(n^{-1}\log n)^{1/2}$ for all $j = 1, \ldots, j_0$ and all $k = j_0 + 1, \ldots, p$, converges to 1 as $n \to \infty$. From this result and the first part of (4.9), it follows that for $1 \le j \le j_0$

$$\hat{r}_j = 1 + \sum_{k\,:\,k \ne j} I(\hat{\theta}_j \le \hat{\theta}_k) = 1 + \sum_{k\,:\,k \le j_0, k \ne j} I(W_j \sigma_j \le W_k \sigma_k) + \Delta_j,$$

where the random variables $W_1, \ldots, W_{j_0}$ are asymptotically independent and distributed as normal $N(0, 1)$, and $P(\Delta_j = 0) \to 1$ as $n \to \infty$.

(ii) In the bootstrap case, it follows from the second formula in (4.3) that

$$(5.3) \quad \hat{r}_j^* = 1 + \sum_{k\,:\,k \le j_0, k \ne j} I\{m^{1/2}(\hat{\theta}_j^* - \hat{\theta}_j) + \Delta_{jk} \le m^{1/2}(\hat{\theta}_k^* - \hat{\theta}_k)\} + \Delta_j^*,$$

where, if $n$ is so large that $\inf_{1 \le j \le j_0} \inf_{j_0 < k \le p}(\theta_j - \theta_k) > 4C_2(m^{-1}\log m)^{1/2}$, then

$$(5.4) \quad \sup_{1 \le k \le j_0} |\Delta_{jk}| \le 2m^{1/2}\Big(\sup_{1 \le j \le j_0} |\hat{\theta}_j - \theta_j| + \sup_{1 \le j_1, j_2 \le j_0} |\theta_{j_1} - \theta_{j_2}|\Big) \to 0,$$

$$(5.5) \quad \begin{aligned} P(\Delta_j^* \ne 0) &\le p \sup_{1 \le k \le p} [P\{|\hat{\theta}_k - \theta_k| > C_2(m^{-1}\log m)^{1/2}\} \\ &\quad + P\{|\hat{\theta}_k^* - \hat{\theta}_k| > C_2(m^{-1}\log m)^{1/2}\}] \to 0. \end{aligned}$$

The convergence in (5.4) is in probability and is a consequence of (4.9), (4.10) and the fact that $m/n \to 0$, and (5.5) follows from (4.7) and the second parts of (4.9) and (4.12). Part (ii) of Theorem 4.1 follows from (5.3)–(5.5).

(iii) Note that

$$E(\hat{r}_j^* | \mathcal{X}, \mathcal{X}_j^*) - 1 = \sum_{k\,:\,k \ne j} P(\hat{\theta}_j^* \le \hat{\theta}_k^* | \mathcal{X}, \mathcal{X}_j^*) = S_1^* + (S_2 + S_3 + S_4^*)\Omega,$$

where $P(0 \le \Omega \le 1) = 1$,

$$S_1^* = \sum_{k=1}^{j_0} P(\hat{\theta}_j^* \le \hat{\theta}_k^* | \mathcal{X}, \mathcal{X}_j^*),$$

$$S_2 = \sum_{k=j_0+1}^{\infty} I\{\theta_j - \theta_k \le 4C_2(m^{-1}\log m)^{1/2}\},$$

$$S_3 = \sum_{k=1}^{p} I\{|\hat{\theta}_k - \theta_k| > C_2(m^{-1}\log m)^{1/2}\},$$

$$S_4^* = \sum_{k=1}^{p} P\{|\hat{\theta}_k^* - \hat{\theta}_k| > C_2(m^{-1}\log m)^{1/2} | \mathcal{X}, \mathcal{X}_j^*\}.$$

In view of (4.10) and (4.11), $S_2 = 0$ for all sufficiently large $n$; by (4.7) and the second part of (4.9), $E(S_3) = o(1)$; and by (4.7) and the second part of (4.12), $E(S_4^*) = o(1)$. Therefore, $E(S_2 + S_3 + S_4^*) = o(1)$; call this result (R). Since, using the independent-component bootstrap, $\mathcal{X}_j^*$ and $\mathcal{X}_k^*$ (for $k \neq j$) are independent conditional on $\mathcal{X}$; and since

$$P(\hat{\theta}_j^* \leq \hat{\theta}_k^* | \mathcal{X}, \mathcal{X}_j^*) = P\{m^{1/2}(\hat{\theta}_j^* - \hat{\theta}_j) + m^{1/2}(\hat{\theta}_j - \hat{\theta}_k) \leq m^{1/2}(\hat{\theta}_k^* - \hat{\theta}_k) | \mathcal{X}\};$$

then it follows from (4.6), the first parts of (4.9) and (4.12), and Kolmogorov's extension theorem, that the joint distribution function of $P(\hat{\theta}_j^* \leq \hat{\theta}_k^* | \mathcal{X}, \mathcal{X}_j^*)$, for $1 \leq k \leq j_0$ and $k \neq j$ (and conditional on $\mathcal{X}$), minus the joint distribution function of $P(Z_j \sigma_j \leq Z_k \sigma_k | Z_j)$ for $1 \leq k \leq j_0$ and $k \neq j$ (for independent standard normal random variables $Z_k$ defined on an enlarged probability space), converges to zero in probability in any integral metric on a compact set. Therefore, the distribution function of $S_1^* + 1$, conditional on $\mathcal{X}$, minus the distribution of $E(R_j | Z_j)$, converges in probability to zero. [Here, $R_j$ is the function of $Z_1, \ldots, Z_{j_0}$ defined at (4.13), and the construction of $Z_1, \ldots, Z_{j_0}$ involves them being measurable in the sigma-field generated by $\mathcal{X} \cup \mathcal{X}_j^*$.] This property, and result (R), together imply (4.14).

To derive the final portion of part (iii) of Theorem 4.1, note that the argument leading to (4.14) implies that

$$E\left\{ \sum_{k=j_0+1}^{\infty} P(\hat{\theta}_j^* \leq \hat{\theta}_k^* | \mathcal{X}, \mathcal{X}_j^*) \right\} = o(1).$$

Therefore,

$$E(\hat{r}_j^{*2} | \mathcal{X}) = \sum_{k_1, k_2 : k_1, k_2 \neq j} E\{P(\hat{\theta}_j^* \leq \hat{\theta}_{k_1}^* | \mathcal{X}, \mathcal{X}_j^*) P(\hat{\theta}_j^* \leq \hat{\theta}_{k_2}^* | \mathcal{X}, \mathcal{X}_j^*) | \mathcal{X}\}$$

$$= \sum_{k_1, k_2 : k_1, k_2 \neq j, 1 \leq k_1, k_2 \leq j_0} E\{P(\hat{\theta}_j^* \leq \hat{\theta}_{k_1}^* | \mathcal{X}, \mathcal{X}_j^*)$$
$$\times P(\hat{\theta}_j^* \leq \hat{\theta}_{k_2}^* | \mathcal{X}, \mathcal{X}_j^*) | \mathcal{X}\} + o_p(1)$$

$$= T_2 + o_p(1),$$

where, for $\ell = 1, 2$,

$$T_\ell = E\left[ \left\{ \sum_{k : k \neq j, 1 \leq k \leq j_0} P(\hat{\theta}_j^* \leq \hat{\theta}_k^* | \mathcal{X}, \mathcal{X}_j^*) \right\}^{\ell} \Big| \mathcal{X} \right] + o_p(1).$$

More simply, $E(\hat{r}_j^* | \mathcal{X}) = T_1 + o_p(1)$. The argument in the previous paragraph can be used to show that $T_1$ and $T_2$ converge in probability to $E\{E(R_j - 1 | Z_j)^2\}$ and $E(R_j - 1)$, respectively. Since $E\{E(R_j - 1 | Z_j)^2\} = E\{(R_j - 1)^2\}$, then $\operatorname{var}(\hat{r}_j^* | \mathcal{X})$ converges in probability to $\operatorname{var}(R_j)$, as had to be proved.

DEPARTMENT OF MATHEMATICS
   AND STATISTICS
UNIVERSITY OF MELBOURNE
MELBOURNE, VIC 3010
AUSTRALIA
E-MAIL: halpstat@ms.unimelb.edu.au
        h.miller@ms.unimelb.edu.au